\documentclass[10pt]{article}
\usepackage{amsmath}
\usepackage{latexsym}
\usepackage{amssymb}
\usepackage{amsfonts}
\usepackage{amsthm}
\title{PARAMETRIZATION OF COSSERAT EQUATIONS}
\author{J.F. Pommaret \\ CERMICS, Ecole Nationale des Ponts et
  Chauss\'ees,\\6/8 Av. Blaise Pascal, 77455 Marne-la-Vall\'ee Cedex 02, France \\          
  e-mail: jean-francois.pommaret@wanadoo.fr \\ 
  home page: http://cermics.enpc.fr/$\sim$pommaret/home.html  }
\date{  } 
\begin{document}
\maketitle

1) {\bf INTRODUCTION}:\\

As a matter of fact, the solution space of many systems of ordinary differential (OD) or partial differential (PD) equations in engineering or mathematical physics "{\it can}/{\it cannot}" be parametrized by a certain number of arbitrary functions behaving like "{\it potentials}". In view of the explicit examples to be met later on, it must be noticed that the parametrizing operator, though often of the first order, may be, on the contrary, of arbitrarily large order. Among the well known examples, we recall that a classical OD {\it control system} is parametrizable if and only if it is controllable (Kalman test of 1969 in [9]). Among PD systems, the electromagnetic (EM) {\it field}, solution of the first set of 4 Maxwell equations, admits a well known first order parametrization by means of the 4-{\it potential} while the EM {\it induction}, solution of the second set of 4 Maxwell equations (in vacuum), also admits a first order parametrization by means of the so-called 4-{\it pseudopotential}. On the contrary, it is now known that, contrary to the EM situation, the set of 10 second order linearized Einstein equations (in vacuum) cannot be parametrized and {\it cannot} therefore be considered as field equations (see [18,24] for more details; see also [30] and http://wwwb.math.rwth-aachen.de/OreModules for a computer  algebra solution). One among the best interesting and useful cases is concerned with continuum mechanics where the first order stress equations (in vacuum) admits a rather simple second order parametrization by means of the single Airy function in dimension 2 and it is not so well known (!) that a much more complicated second order parametrization can be achieved in dimension $n\geq 2$ by means of $n^2(n^2-1)/12$ arbitrary functions.\\

It is now known that all the above problems are only very particular cases of a {\it much more sophisticated} (!) and general situation involving the {\it formal theory} of systems of PD equations pioneered by D.C. Spencer in 1970 [26] (jet theory, diagram chasing, differential sequences,...)(see [22,23] for more details) and {\it differential modules} in the framework of "{\it algebraic analysis}" pioneered about at the same time by V.P. Palamodov [17] for the constant coefficients case, M. Kashiwara [10] and B. Malgrange [15] for the variable coefficients case (see [23,24] for more details and also consult Zentralblatt Zbl 1079.93001 for a review). The corresponding {\it duality theory,} that is at the heart of all the previous examples and will be a central tool in this paper, highly depends on (hard) {\it homological algebra techniques} (localization, resolutions, extension modules,...) {\it which cannot be avoided}. \\

The purpose of this paper is to apply these techniques in a way as simple and self-contained as possible in order to give a positive and explicit answer concerning the possibility to exhibit a {\it first order} parametrization of the stress/couple-stress equations met in the study of  Cosserat media. At the same time, as a corollary of the homological test, we shall give {\it for the first time} the reason for which the compatibility conditions (CC) for the deformation tensor in classical elasticity theory are {\it second order} while the corresponding CC for Cosserat fields [11] are only {\it first order} and explain why this order is equal to the order of the corresponding parametrization. \\

At the end of the paper, we shall give hints in order to explain why, though the "{\it fields}" and their CC in classical and Cosserat elasticity theories look like {\it completely different at first sight}, therefore providing different presentations of the corresponding field equations, nevertheless the possibility to obtain a parametrization in one framework {\it necessarily} implies the possibility to have a parametrization in the other framework and vice-versa. Though striking it may look like in such an engineering background, this {\it totally not evident result}, which was not known up to now, is one of the simplest consequences of a delicate theorem of homological algebra (see [25], Th 7.22, p 212). In particular {\it the reader must look at the first section below with care, even though it does not seem to have anything to do with Cosserat media}.\\

\noindent
2)  {\bf MOTIVATIONS}:\\

     In the middle of the last century, {\it commutative algebra}, namely the study of modules over rings, was facing a very subtle problem, the resolution of which led to the modern but difficult {\it homological algebra} with sequences and diagrams. Roughly, one can say that the problem was essentially to study properties of finitely generated modules not depending on the "{\it presentation}" of these modules by means of generators and relations. This {\it very hard step}" is based on homological/cohomological methods like the so-called "{\it extension}" modules which cannot therefore be avoided [1,7,16,25].\\

\noindent
{\bf EXAMPLE 2.1}: In order to sketch this problem, let us present two simple examples. In the first case with standard notations, everybody will understand at once that integrating the second order OD equation $\ddot{y}=0$ with $m=n=1$ is equivalent to integrating the system of two first order OD equations ${\dot{y}}^1-y^2=0, {\dot{y}}^2=0$. However, even with $m=n=2$ and the same two unknowns $u,v$ in both cases, it is not evident at all that integrating the second order PD equation $d_{12}u-d_{22}v-u=0$ is equivalent to integrating the system of two fourth order PD equations $d_{1122}u-d_{1222}v-d_{22}v-u=0, d_{1112}u-d_{1122}v-d_{11}u=0$ (exercise !).\\
     
     As before, using now rings of "{\it differential operators}" instead of polynomial rings led to {\it differential modules} and to the challenge of adding the word "{\it differential}" in front of concepts of commutative algebra. Accordingly, not only one needs properties not depending on the presentation as we just explained but also properties not depending on the coordinate system as it becomes clear from any application to mathematical or engineering physics where tensors and exterior forms are always to be met like in the space-time formulation of electromagnetism. Unhappily, no one of the previous techniques for OD or PD equations could work !.\\
   
     By chance, the intrinsic study of systems of OD or PD equations has been pioneered in a totally independent way by D. C. Spencer and collaborators after 1960 [26], using jet theory and diagram chasing in order to relate differential properties of the equations to algebraic properties of their symbol, a technique superseding the "{\it leading term}" approach of Janet in 1920 [8] or Gr\"{o}bner in 1940 [6] but quite poorly known by the mathematical community, even today.\\
     
     Accordingly, it was another challenge to unify the "{\it purely differential}" approach of Spencer with the "{\it purely algebraic}" approach of commutative algebra, having in mind the necessity to use the previous homological algebraic results in this new framework. This sophisticated mixture of differential geometry and homological algebra, now called "{\it algebraic analysis}", has been achieved after 1970 by V. P. Palamodov for the constant coefficient case [17], then by M. Kashiwara [10] and B. Malgrange [15] for the variable coefficient case.\\

Let $k$ be field of characteristic zero and $\chi=({\chi}_1,...,{\chi}_n)$ be indeterminates over $k$. We introduce the ring $A=k[{\chi}_1,...,{\chi}_n]$ of polynomials with coefficients in $k$. When a given system of linear PD equations of order $q$ is given, it defines by residue a differential module $M$ over the underlying ring $D$ of differential operators. More precisely, let us introduce $n$ commuting derivatives $d_1,...,d_n$ for which $k$ should be a field of constants and define the ring $D=k[d]=k[d_1,...,d_n]$ of differential operators with coefficients in $k$. Then $D$ and $A$ are isomorphic by $d_i\leftrightarrow {\chi}_i$. However, the (non-commutative) situation for a differential field $K$ with subfield of constants $k$ escapes from the previous (commutative) approach and must be treated "{\it by its own}". For this, let $\mu=({\mu}_1,...,{\mu}_n)$ be a multi-index with {\it length} $\mid\mu\mid={\mu}_1+ ... +{\mu}_n$. We set $\mu +1_i=({\mu}_1,...,{\mu}_{i-1},{\mu}_i+1,{\mu}_{i+1},...,{\mu}_n)$ and we say that $\mu$ is of {\it class} $i$ if ${\mu}_1=...={\mu}_{i-1}=0, {\mu}_i\neq 0$. Accordingly, any operator $P=a^{\mu}d_{\mu}\in D$ acts on the unknowns $y^k$ for $k=1,...,m$ as we may set $d_{\mu}y^k=y^k_{\mu}$ with $y^k_0=y^k$ and introduce the {\it jet coordinates} $y_q=\{y^k_{\mu}\mid k=1,...,m; 0\leq \mid \mu \mid \leq q\}$. It follows that, if a system of PD equations can be written in the form ${\Phi}^{\tau}\equiv a^{\tau\mu}_ky^k_{\mu}=0$ with $a\in K $, we may introduce the filtred differential module $M=Dy/D\Phi$. Then we define the (formal) {\it prolongation} of ${\Phi}^{\tau}$ with respect to $d_i$ to be $d_i{\Phi}^{\tau}\equiv a^{\tau\mu}_ky^k_{\mu +1_i}+ {\partial}_ia^{\tau\mu}_ky^k_{\mu}$ and induce maps $d_i:M\rightarrow M: {\bar{y}}^k_{\mu}\rightarrow {\bar{y}}^k_{\mu +1_i}$ by residue.\\

First of all, we sketch the technique of "{\it localization}" in the case of OD equations, comparing to the situation met in classical control theory where $n=1$. For this, setting as usual $d=d_1=d/dt=dot$, we may introduce (formal) unknowns $y^1,...,y^m$ and set $Dy=Dy^1+...+Dy^m\simeq D^m$. If we have a given system $\Phi=0$ of OD equations of order $q$, a basic question in control theory is to decide whether the control system is "{\it controllable}" or not. It is not our purpose to discuss here about such a question (see [23,24] for more details) but we just want to state the final formal test in terms of a property of the differential module $M=Dy/D\Phi$. Care must be taken that in the sequel, for simplicity and unless needed, we shall not always put a "{\it bar}" on the residual image of $y$ in the canonical projection $Dy\rightarrow M$. We explain our goal on an example.\\

\noindent
{\bf EXAMPLE 2.2}: With $m=3$ and a constant parameter $a$, we consider the first order system ${\Phi}^1\equiv {\dot{y}}^1-ay^2-{\dot{y}}^3=0, {\Phi}^2\equiv y^1-{\dot{y}}^2+{\dot{y}}^3=0$. Any engineer should want to apply Laplace transform $\hat{y}(s)={\int}^{\infty}_0e^{st}y(t)dt$ to this system. However, using the integration by part formula ${\int}^{\infty}_0e^{st}{\dot{y}}(t)dt=[e^{st}y(t)]^{\infty}_0-s\hat{y}(s)$ we should eventually need to know $y(0)$ though the Kalman test of controllability is purely formal as it only deals with ranks of matrices [9]. Since a long time we had in mind that setting $y(0)=0$ was not the right way and that Laplace transform could be superseded by another purely formal technique. For this, let us replace "{\it formally}" $d$ by the purely algebraic symbol $\chi$ whenever it appears and obtain the system of {\it linear equations} :\\
\[ \chi y^1-ay^2-\chi y^3=0, 1y^1-\chi y^2+\chi y^3=0 \Rightarrow y^1=\frac{\chi (\chi +a)}{{\chi}^2-a} y^3, y^2=\frac{\chi (\chi +1)}{{\chi}^2-a} y^3 \]
but we could have adopted a different choice for the only arbitrary unknown. At this step there are only two possibilities :\\
$\bullet a\neq 0,1 \Rightarrow $no "{\it simplification}" may occur and, getting rid of the common denominator, we get an algebraic parametrization leading to a differential parametrization as follows:\\
\[ y^1=\chi (\chi +a) z, y^2=\chi (\chi +1) z, y^3=({\chi}^2-a)z \Rightarrow y^1=\ddot{z}+a\dot{z}, y^2=\ddot{z}+\dot{z}, y^3=\ddot{z}-az  \]
$\bullet a=0$ or $a=1 \Rightarrow$  a "{\it simplification}" may occur and no parametrization can be found. For example, with $a=0$, setting $z=y^1-y^3$ we get $\chi z=0$ that is to say $\dot{z}=0$.\\

Recapitulating, we discover that a control system is controllable if and only if one cannot get any {\it autonomous} element satisfying an OD equation {\it by itself}. For understanding such a result in an algebraic manner, let $M$ be a module over an integral domain $A$ containing 1. A subset $S\subset A$ is called a {\it multiplicative subset} if $1\in S$ and $\forall s,t \in S \Rightarrow st \in S$. Moreover, we shall need/use the {\it Ore condition} on $S$ and $A$, namely $aS\cap sA\neq \emptyset , \forall a\in A, s\in S$.\\

\noindent
{\bf DEFINITION 2.3}: $S^{-1}A=\{s^{-1}a{\mid}s\in S,a\in A/\sim\}$ with $s^{-1}a\sim t^{-1}b \Leftrightarrow \exists u,v\in A, us=vt\in S , ua=vb$.\\

Next, for any module $M$ over $A$, we define $S^{-1}M=S^{-1}A{\otimes}_AM$ and $t_S(M)=\{x\in M\mid \exists s\in S, sx=0\} $ in the exact sequence $0\rightarrow t_S(M)\rightarrow M \rightarrow S^{-1}M$ 
where  the last morphism is $x \rightarrow 1^{-1}x$.\\

\noindent
{\bf EXAMPLE 2.4}: $S=A-\{ 0 \} \Rightarrow S^{-1}A=Q(A) $ field of fractions of $A$ and we introduce the {\it torsion submodule} $t_S(M)=t(M)=\{ x\in M \mid \exists 0\neq a\in A, ax=0\}$ of $M$. \\

In the case of a torsion-free module, reducing to the same denominator as in the control example, we obtain the following classical proposition amounting to exhibit a parametrization [23,25]. However, the reader must notice that it is unusefull in actual practice as {\it one needs a test} (like the Kalman test) for checking the torsion-free condition. {\it This will be the hard part of the job} !.\\

\noindent
{\bf PROPOSITION 2.5}: When $M$ is finitely generated and $t(M)=0$, from the inclusion $M\subset Q(A){\otimes}_A M$, we deduce that there exists a finitely generated free module $F$ with $M\subset F$.\\ 

\noindent
{\bf EXAMPLE 2.6}: In 2-dimensional classical elasticity, let us consider the well known stress equations:\\
\[  {\partial}_1{\sigma}^{11}+{\partial}_2{\sigma}^{21}=0, {\partial}_1{\sigma}^{12}+{\partial}_2{\sigma}^{22}=0\]
where we must keep in mind that the stress tensor density is symmetric, that is ${\sigma}^{12}={\sigma}^{21}$. Replacing ${\partial}_i$ by ${\chi}_i$, we may localize and obtain:\\
\[    {\chi}_1{\sigma}^{11}+{\chi}_2{\sigma}^{21}=0, {\chi}_1{\sigma}^{12}+{\chi}_2{\sigma}^{22}=0  \]
Reducing the fractions to the same denominator, we get:\\
\[  {\sigma}^{11}=-\frac{{\chi}_2}{{\chi}_1}{\sigma}^{21}=-\frac{({\chi}_2)^2}{{\chi}_1{\chi}_2}{\sigma}^{12}, {\sigma}^{22}=-\frac{{\chi}_1}{{\chi}_2}{\sigma}^{12}=-\frac{({\chi}_1)^2}{{\chi}_1{\chi}_2} {\sigma}^{12} \]
and obtain therefore the subvector space over $\mathbb{Q}({\chi}_1,{\chi}_2)$:\\
\[ {\sigma}^{11}=({\chi}_2)^2\phi, {\sigma}^{12}={\sigma}^{21}=-{\chi}_1{\chi}_2\phi, {\sigma}^{22}=({\chi}_1)^2\phi \]
a result providing at once the well known parametrization by the Airy function:\\
\[{\sigma}^{11}={\partial}_{22}\phi, {\sigma}^{12}={\sigma}^{21}=-{\partial}_{12}\phi, {\sigma}_{22}={\partial}_{11}\phi  \]
It may be interesting to compare this purely formal approach to the standard analytic aproach presented in any textbook along the following way. From the first stress equation and Stokes identity for the curl, there exists a function $\varphi$ such that ${\sigma}^{11}={\partial}_2\varphi, {\sigma}^{21}=-{\partial}_1\varphi$. Similarly, from the second stress equation, there exists a function $\psi$ such that ${\sigma}^{22}={\partial}_1\psi, {\sigma}^{12}=-{\partial}_2\psi$. Finally, from the symmetry of the stress, there exists a function $\phi$ such that $\varphi={\partial}_2\phi, \psi={\partial}_1\phi$ and we find back the same parametrization of course. The reader must notice that, in this example, one can check that the parametrization does work but no geometric inside can be achieved in arbitrary dimension $n\geq 2$ even though exactly the same procedure can be applied.\\

\noindent
{\bf EXERCISE 2.7}: Apply similarly the localization technique in the case of the two sets of Maxwell equations and compare to the standard analytic approach to be found in any textbook.\\

The extension of the above results to the non-commutative case $D=K[d]$ where $K$ is a differential field with $n$ commuting derivations ${\partial}_1,...,{\partial}_n$  can be achieved but is more delicate 
[  23].\\

\noindent
{\bf EXERCISE 2.8}: When $a=a(t)$ in Example 1.6, we let the reader prove that the controllability condition is now the Ricatti inequality $\dot{a} + a^2-a\neq 0$ in a coherent way with the constant coefficient case already considered.\\

\noindent
Taking into account the works of Janet and Spencer, the study of systems of PD equations cannot be achieved without understanding {\it involution} and we now explain this concept. For this, changing linearly the derivations if necessary, we may successively solve the maximum number of equations with respect to the jets of order $q$ and class $n$, class $(n-1)$,..., class 1. At each order, a certain number of jets called {\it principal} can therefore be expressed by means of the other jets called {\it parametric}. Moreover, for each equation of order $q$ and class i, $d_1,...,d_i$ are called {\it multiplicative} while $d_{i+1},...,d_n$ are called {\it nonmultiplicative} and $d_1,...,d_n$ are nonmultiplicative for all the remaining equations of order $\leq q-1$.\\

\noindent
{\bf DEFINITION 2.9}: The system is said to be {\it involutive} if each prolongation with respect to a nonmultiplicative derivation is a linear combination of prolongations with respect to the multiplicative ones. Using Spencer cohomology, one can prove that such a definition is in fact intrinsic [22,23,26].\\

\noindent
{\bf EXAMPLE 2.10}: The system $y_{11}=0, y_{13}-y_2=0$ is not involutive. Effecting the permutation $(1,2,3)\rightarrow (3,2,1)$, we get the new system $y_{33}=0, y_{13}-y_2=0$. As $d_1y_{33}-d_3(y_{13}-y_2)=y_{23}$ and $d_1y_{23}-d_2(y_{13}-y_2)=y_{22}$, the new system $y_{33}=0, y_{23}=0, y_{22}=0, y_{13}-y_2=0$ is involutive with 1 equation of class 3, 2 equations of class 2 and 1 equation of class 1.\\

\noindent
{\bf EXAMPLE 2.11}: The Killing system $y^j_i+y^i_j=0$ is {\it not} involutive but the first prolongation $y^j_i+y^i_j=0, y^k_{ij}=0$ is involutive. This is the reason for which the Riemann tensor is a first order expression in the metric and Christoffel symbols and thus second order in the metric alone (for more details, see [22], p 249-258).\\

\noindent
{\bf APPLICATION 2.12}: $t(M)=M$ if and only if the number of equations of class $n$ is $m$. Otherwise there is a strict inclusion $t(M)\subset M$.\\

\noindent
{\bf PROPOSITION 2.13}: ([22,26]) The following recipe (already used in the example at the beginning of this section) will allow to bring an involutive system of order $q$ to an equivalent (isomorphic modules) involutive system of order 1 {\it with no zero order equations} called 
{\it Spencer form}:\\
1) Use all parametric jets up to order $q$ as new unknowns.\\
2) Make one prolongation.\\
3) Substitute the new unknowns.\\

\noindent
3)  {\bf GROUP FOUNDATION OF MECHANICS}:\\

   This section, which is a summary of results already obtained in [21,22], is provided for fixing the notations and the techniques leading to various (linearized) {\it differential sequences}. All the results presented are local ones. A corresponding non-linear does exist but is out of the scope of this paper [12,21].\\
Let $X$ be a manifold of dimension $n$ with local coordinates $x=(x^1,...,x^n)$ and latin indices $i,j=1,...,n$. We denote by $T=T(X)$ the tangent bundle to $X$ and by $T^*=T^*(X)$ the cotangent bundle to $X$ while ${\wedge}^rT^*$ is the bundle of $r$-forms on $X$. Also, we denote by $J_q(T)$ the $q$-{\it jet bundle} of $T$, that is to say the vector bundle over $X$ having the same transition rules as a vector field and its derivatives up to order $q$ under any change of local coordinates on $X$. Let now $G$ be a Lie group of dimension $p$ with identity $e$, local coordinates $a=(a^1,...,a^p)$ and greek index $\tau$. We denote by ${\cal{G}}=T_{e}(G)$ the corresponding Lie algebra with vectors denoted by the greek letters $\lambda$. As usual, we shall identify a map $a:X\rightarrow G$ called {\it gauging} of $G$ over $X$, with its graph $X\rightarrow X\times G$ which is a section of a trivial principal bundle, and, similarly, use the same notation for a bundle and its sheaf of (local) sections as the background will always tell the right choice. In particular, when differential operators are involved, the sectional point of view must automatically be used. Such a convention allows to greatly simplify the notations at the expense of a slight abuse of language.\\

\noindent
{\bf DEFINITION 3.1}: A {\it Lie group of transformations} of a manifold $X$ is a lie group $G$ with an {\it action} of $G$ on $X$ better defined by its graph $X\times G \rightarrow X\times X: (x,a)\rightarrow (x,y=ax=f(x,a))$ with the properties that $a(bx)=(ab)x$ and $ex=x, \forall x\in X, \forall a,b\in G$. \\

    It is sometimes useful to distinguish the {\it source} $x$ from the {\it target} $y$ by introducing a copy $Y$ of $X$ with local coordinates $y=(y^1,...,y^n)$. Such groups of transformations have first been studied by S. Lie in 1880. Among basic examples when $n=1$ we may quote the {\it affine group} $y=ax+b$ and the {\it projective group} $y=(ax+b)/(cx+d)$ of transformations of the real line. When $n=3$ we may quote the {\it group of rigid motions} $y=ax+b$ where now $a$ is an orthogonal $3\times 3$ matrix and $b$ is a vector. Such a group is known to preserve the {\it euclidean metric} $\omega=({\omega}^{ij}={\omega}^{ji})$ and thus the quadratic form $ds^2=(dx^1)^2+(dx^2)^2+(dx^3)^2={\omega}^{ij}dx^idx^j$. When $n=4$ we may quote the {\it conformal group} of space-time with $15$ parameters (4 translations, 6 rotations, 1 dilatation, 4 elations) preserving the Minkowski metric $\omega$ or the quadratic form $ds^2=(dx^1)^2+(dx^2)^2+(dx^3)^2-c^2(dt)^2$ up to a function factor, where now $c$ is the speed of light and $t$ the time. Among the subgroups, we may consider the {\it Weyl group} with $11$ parameters preserving $\omega$ up to a constant factor and the {\it Poincar\'e group} with $10$ parameters preserving $\omega$.\\

Only ten years later, in 1890, S. Lie discovered that the Lie groups of transformations were only examples of  a wider class of groups of transformations, first called {\it infinite groups} but now called {\it Lie pseudogroups}.\\

\noindent
{\bf DEFINITION 3.2}: A Lie pseudogroup $\Gamma$ of transformations of a manifold $X$ is a group of transformations $y=f(x)$ solutions of a (in general nonlinear) system of PD equations, also called system of {\it finite Lie equations}.

Setting now $y=x+t{\xi}(x)+...$ and passing to the limit for $t\rightarrow 0$, that is to say linearizing the defining system of finite Lie equations around the q-jet of the identity $y=x$, we get a linear system $R_q\subset J_q(T)$ for vector fields, also called system of {\it infinitesemal Lie equations}, with solutions $\Theta\subset T$ satisfying $[\Theta,\Theta]\subset \Theta$. It can be proved, for the same testing type reasons, that such a system may be endowed with a Lie algebra bracket {\it on sections} ${\xi}_q:(x)\rightarrow (x,{\xi}^k(x), {\xi}^k_i(x),{\xi}^k_{ij}(x),...)$ that we shall quickly define (see [22] for more details and compare to [12]). Such a bracket on sections transforms $R_q$ into a {\it Lie algebroid} and we have $[R_q,R_q]\subset R_q$. Let us first define by bilinearity $\{j_{q+1}(\xi),j_{q+1}(\eta)\}=j_q([\xi,\eta]), \forall \xi,\eta\in T$. Introducing the {\it Spencer operator} $D:R_{q+1}\rightarrow T^*\otimes R_q:{\xi}_{q+1}\rightarrow j_1({\xi}_q)-{\xi}_{q+1}$ with local components $({\partial}_i{\xi}^k-{\xi}^k_i,{\partial}_i{\xi}^k_j-{\xi}^k_{ij},...)$, we obtain the following general formula at order $q$:\\
\[   [{\xi}_q,{\eta}_q]=\{{\xi}_{q+1},{\eta}_{q+1}\}+i(\xi)D{\eta}_{q+1}-i(\eta)D{\xi}_{q+1}  ,\forall {\xi}_q,{\eta}_q\in R_q\]
where $i( )$ is the interior multiplication of a 1-form by a vector, and we let the reader check that such a definition no longer depends on the "lifts" ${\xi}_{q+1},{\eta}_{q+1}$ over ${\xi}_q,{\eta}_q$.\\

\noindent
{\bf EXAMPLE 3.3}: (Affine transformations)  $n=1,q=2, X={\mathbb{R}}^3$\\
With evident notations, the system of finite Lie equations is defined by the single second order linear OD equation $y_{xx}=0$. Similarly, the sections of $R_2$ are defined by ${\xi}_{xx}(x)=0$. Accordingly, the components of $[{\xi}_2,{\eta}_2]$ at order zero, one and two are defined by the totally unusual successive formulas:\\
\[    [\xi,\eta]=\xi{\partial}_x\eta-\eta{\partial}_x\xi     \]
\[    ([{\xi}_1,{\eta}_1])_x=\xi{\partial}_x{\eta}_x-\eta{\partial}_x{\xi}_x    \]
\[    ([{\xi}_2,{\eta}_2])_{xx}={\xi}_x{\eta}_{xx}-{\eta}_x{\xi}_{xx}+\xi{\partial}_x{\eta}_{xx}-\eta{\partial}_x{\xi}_{xx}   \]
It follows that ${\xi}_{xx}=0,{\eta}_{xx}=0\Rightarrow ([{\xi}_2,{\eta}_2])_{xx}=0$ and thus $[R_2,R_2]\subset R_2$.\\

\noindent
{\bf EXAMPLE 3.4}: (Projective transformations)  $n=1, q=3, X={\mathbb{R}}^3$\\
The system of finite Lie equations is defined by the single third order nonlinear OD equation $(y_{xxx}/y_x)-\frac{3}{2}{(y_{xx}/y_x)}^2=0$ and the sections of $R_3$ are defined by ${\xi}_{xxx}(x)=0$. The formulas for the bracket of Lie algebroid $[R_3,R_3]\subset R_3$ can be derived similarly but involve many more terms.\\

\noindent
{\bf EXAMPLE 3.5}: (Volume preserving transformations)  $n$ arbitrary, $q=1, X={\mathbb{R}}^n$\\
The only non-linear finite Lie equation is ${\partial}(y^1,...,y^n)/{\partial}(x^1,...,x^n)=det(y^k_i)=1$ and the sections of $R_1$ are defined by the single relation ${\xi}^i_i=0$. Accordingly, we obtain:\\
\[   ([{\xi}_1,{\eta}_1])^k_i={\xi}^r_i{\eta}^k_r-{\eta}^r_i{\xi}^k_r+{\xi}^r{\partial}_r{\eta}^k_i-{\eta}^r{\partial}_r{\xi}^k_i       \]
When summing on $k$ and $i$, the first two terms disappear (as in Example 3.3 !) and we get therefore $[R_1,R_1]\subset R_1$. We invite the reader to compare this result with the usual way on solutions where one defines $\Theta$ as the kernel of the Lie derivative ${\cal{L}}(\xi)\omega$ of the (volume) n-form $\omega=dx^1\wedge ...\wedge dx^n$ with respect to $\xi$ and then uses the well known formula $[{\cal{L}}(\xi),{\cal{L}}(\eta)]\omega={\cal{L}}([\xi,\eta])\omega , \forall \xi,\eta\in T$ in order to obtain $[\Theta,\Theta]\subset \Theta$.\\

    Introducing a basis of ${\wedge}^rT^*$ made by the $dx^I=dx^{i_1}\wedge ... \wedge dx^{i_r}$ with $I=(i_1<...<i_r)$, we may define the {\it exterior derivative} $d:{\wedge}^rT^*\rightarrow {\wedge}^{r+1}T^*$ by setting $\omega={\omega}_Idx^I  \rightarrow d\omega ={\partial}_i{\omega}_Idx^i\wedge dx^I$ and one easily checks $d^2=d\circ d=0$. This is sufficient in order to define the (canonical linear) {\it gauge sequence}:\\
\[   {\wedge}^0T^*\otimes {\cal{G}}\stackrel{d}{\rightarrow}{\wedge}^1T^*\otimes {\cal{G}}\stackrel{d}{\rightarrow}{\wedge}^2T^*\otimes {\cal{G} }     \]
which is just the tensor product by $\cal{G}$ of a part of the Poincar\'e sequence for the exterior derivative.\\

    However, there are other differential sequences to be found in the literature and that we did not speak about, namely the {\it Janet sequence}, which is for sure the best known differential sequence, and the {\it Spencer sequence}. For short, if $E,F,...$ denote vector bundles over $X$, we use the same letters for the corresponding sets (sheaves to be exact) of sections and such an interpretation must be used whenever operators are involved. Starting from a vector bundle $E$ (for example $T$) and a linear differential operator ${\cal{D}}:E\rightarrow F:\xi\rightarrow \eta$, if we want to solve the linear system with second member ${\cal{D}}\xi=\eta$ even locally, one needs "{\it compatibility conditions}" (CC) in the form ${\cal{D}}_1\eta=0$. Denoting now $F$ by $F_0$, we may therefore look for an operator ${\cal{D}}_1:F_0\rightarrow F_1:\eta\rightarrow \zeta$ and so on. Under the assumption that $\cal{D}$ is involutive (!), the french mathematician M. Janet has proved in 1920 that such a chain of operators ends after $n$ steps and we obtain the (canonical linear) {\it Janet sequence}, namely [8,22]:\\
 \[   0\rightarrow \Theta\rightarrow E\stackrel{\cal{D}}{\rightarrow}F_0\stackrel{{\cal{D}}_1}{\rightarrow}F_1\stackrel{{\cal{D}}_2}{\rightarrow}...\stackrel{{\cal{D}}_n}{\rightarrow} F_n\rightarrow 0  \]
 where, apart from the first operator which is of order $q$, the $n$ operators following it are first order and  involutive. The (canonical linear) {\it Spencer sequence} is the Janet sequence for the first order Spencer form, namely:\\
 \[  0\rightarrow \Theta \stackrel{j_q}{\rightarrow} C_0\stackrel{D_1}{\rightarrow}C_1\stackrel{D_2}{\rightarrow} ... \stackrel{D_n}{\rightarrow}C_n\rightarrow 0  \]
 where $C_0=R_q$ and the first order involutive operators $D_i$ are induced by the Spencer operator. It follows that we only have at our disposal for any application where group theory seems to be involved, three linear differential sequences, namely the Janet sequence, the Spencer sequence and the gauge sequence. As these sequences are made by quite different operators, {\it the use of one excludes the use of the others}.\\
 
 In order to escape from this dilemna and for the sake of clarifying the key idea of the brothers Cosserat by using these new mathematical tools, we shall explain, in a way as elementary as possible while using only the linear framework, why {\it the Janet sequence and the gauge sequence cannot be used in continuum mechanics}. By this way we hope to convince the reader about the need to use another differential sequence, namely the SPENCER SEQUENCE, though striking it could be. Also we shall use very illuminating examples in order to illustrate our comments.\\
First of all we exhibit the isomoprphism existing between the linear gauge sequence and the linear Spencer sequence. For this, if {\it now} $G$ {\it acts on} $X$ with a basis ${{\xi}_{\tau}=\{\xi}^k_{\tau}{\partial}_k\} $ of infinitesimal generators, we may introduce the map:\\
\[     {\wedge}^0T^*\otimes {\cal{G}}\rightarrow J_q(T): {\lambda}^{\tau}(x) \rightarrow {\lambda}^{\tau}(x){\partial}_{\mu}{\xi}^k_{\tau}(x)     \]
It is known [21 ,p. 308] that this map becomes injective for $q$ large enough and we may call $R_q$ its image for such a $q$. It follows from its definition that $R_q\simeq R_{q+1}$ is a system of infinitesimal Lie equations of finite type and we get for the Spencer operator (care to the notation) [22,26]:\\
\[     D:R_{q+1}\rightarrow T^*\otimes R_q:{\xi}_{q+1}\rightarrow ({\partial}_i{\xi}^k_{\mu}-{\xi}^k_{\mu+1_i}={\partial}_i{\lambda}^{\tau}(x){\partial}_{\mu}{\xi}^k_{\tau}(x))  \]
Accordingly, the linear gauge sequence is isomorphic to the linear Spencer sequence:\\
\[  0\rightarrow \Theta \rightarrow {\wedge}^0T^*\otimes R_q\stackrel{D}{\rightarrow} {\wedge}^1T^*\otimes R_q\stackrel{D}{\rightarrow} {\wedge}^2T^*\otimes R_q  \]
the three isomorphisms being induced by the (local) isomorphism $X\times {\cal{G}}\rightarrow R_q$ just described above. It is essential to notice that, though the linear Spencer sequence and the isomorphisms crucially depend on the action, by a kind of "{\it miracle}" the linear gauge sequence no longer depends on the action (look at the kernel of the first "$d$").\\
This result proves that {\it the linear Spencer sequence generalizes the linear gauge sequence},  with the major gain that it can be used even for Lie pseudogroups of transformations that are not coming from Lie groups of transformations.\\

\noindent
{\bf REMARK 3.6}: When $n=3$ and we deal with the Lie group of rigid motions, though surprising it may look like, the formal adjoint of the (first) Spencer operator is {\it exactly} describing the so-called {\it stress} and {\it couple-stress} equations found by the brothers Cosserat as we shall explain below. Though this result is in perfect agreement with the piezzoelectric or photoelastic coupling of elasticity and electromagnetism [23], it CONTRADICTS gauge theory where the lagrangians are functions on ${\wedge}^2T^*\otimes \cal{G}$ and NOT on $T^*\otimes\cal{G}$ as above [29].\\

  Let us consider a (finite) volume ${\int}_VdV$ in ${\mathbb{R}}^3$ limited by a (closed) surface $S={\int}_SdS$ and let us introduce the outside unit normal (pseudo) vector $\vec{n}=(n_j)$ on $S$. Let us now suppose that the surface element $dS$ is acted on by the outside with a force $d\vec{F}=\vec{\sigma}dS$ and a couple $d\vec{C}=\vec{\mu}dS$, where both $\vec{\sigma}$ and $\vec{\mu}$ linearly depend on $\vec{n}$ through the {\it stress} tensor density $\sigma=({\sigma}^{ij})$ and the {\it couple-stress} tensor density $\mu=({\mu}^{r,ij}=-{\mu}^{r,ji})$. It must be noticed that, using the standard Cauchy tetrahedral device, there is no reason "{\it a-priori}" to suppose that the stress tensor is symmetric. We also suppose that   the volume element $dV$ is acted on by (see later on for the sign) a force $-\vec{f}dV$ and a momentum $-\vec{m}dV$ with $\vec{f}=(f^j)$ and $\vec{m}=(m^{ij}=-m^{ji})$.\\
  Our purpose is now to study the equilibrium of the corresponding torsor fields with respect to an arbitrary cartesian frame $0x^1x^2x^3$.\\
   The equilibrium of forces is satisfied if we have the relation:\\
   \[      {\int}_S\vec{\sigma}dS-{\int}_V\vec{f}dV=0  \Rightarrow {\int}_S{\sigma}^{ij}n_idS-{\int}_Vf^jdV=0\]
Using Stokes formula, this is equivalent to the well known {\it stress equations}:\\
\[        {\partial}_i{\sigma}^{ij}=f^j   \]
This result shows that the surface density of forces $\vec{\sigma}$ is equivalent, from the point of view of force equilibrium, to a volume density of forces $\vec{f}$ and this interpretation explains the sign adopted.\\
   Finally, the equilibrium of forces being satisfied, it is known that the equilibrium of momenta is also satisfied if it is satisfied with respect to an arbitrarily chosen cartesian frame. Hence, introducing the vector $\vec{r}=(x^1,x^2,x^3)$, the equilibrium of momenta is satisfied if we have the relation:\\
   \[     {\int}_S(\vec{\mu}+\vec{r}\wedge \vec{\sigma})dS-{\int}_V(\vec{m}+\vec{r}\wedge \vec{f})dV=0 \]
   Projecting onto the axis $Ox^3$, we obtain:\\
   \[  {\int}_S({\mu}^{r,12}+x^1{\sigma}^{r2}-x^2{\sigma}^{r1})n_rdS-{\int}_V(m^{12}+x^1f^2-x^2f^1)dV=0\]
   Using again Stokes formula and the previous stress equations, we obtain the {\it couple-stress equations}:\\ 
   \[      {\partial}_r{\mu}^{r,ij}+{\sigma}^{ij}-{\sigma}^{ji}=m^{ij}  \]
 This result shows that the surface density of forces $\vec{\sigma}$ and couples $\vec{\mu}$ is equivalent, from the point of view of torsor equilibrium, to a volume density of forces $\vec{f}$ and to a volume density of momenta $\vec{m}$, provided the preceding stress and couple-stress equations are satisfied, and this interpretation explains the sign adopted.\\
The combination of the stress {\it and} couple-stress equations have first been exhibited by E. and F. Cosserat in 1909 [4,5,p137] {\it without any static equilibrium experimental background} and we now explain the key argument leading to the same equations just from group theoretical arguments. Of course, most of the engineering continua such as steel, concrete, glass, wate,... have the specific "{\it constitutive laws}"  $\mu=0, m=0$ and we obtain therefore ${\sigma}^{ij}={\sigma}^{ji}$, that is the stress tensor is symmetric, a situation not always encountered in liquid crystals.\\

First of all, for the reader not familiar with the Spencer operator, we exhibit a similar result in a quite simpler 1-dimensional situation that will allow to recapitulate all the previous results..\\

\noindent
{\bf EXAMPLE 3.7}: Let us consider the Lie group of affine transformations of the real line defined by the group action $y=a^1x+a^2$. The corresponding 2-dimensional Lie group $G$ has coordinates $a=(a^1,a^2)$ and the group composition law is $ab=(a^1,a^2)(b^1,b^2)=(a^1b^1,a^1b^2+a^2)$ with inverse law $a^{-1}=(1/a^1,-a^2/a^1)$. A basis of infinitesimal generators of the action may be obtained with ${\xi}^1=x\frac{\partial}{\partial x}, {\xi}^2=\frac{\partial}{\partial x}$ and we have $[{\xi}^1,{\xi}^2]=-{\xi}^2$. As the exterior derivative is a linear involutive first order homogeneous operator, both with its formal adjoint, it could not have anything to do with the stress and couple-stress equations previously exhibited.\\
{\it Let us now deal with the Spencer sequence instead of the gauge sequence in this framework}.\\
First of all, we may consider the above Lie group of transformations as a Lie pseudogroup defined by the second order system of finite Lie equations $y_{xx}=0$. The corresponding system $R_2\subset J_2(T)$  of infinitesimal Lie equations is ${\xi}_{xx}=0$ and the isomorphism between the gauge sequence and the Spencer sequence is induced by the map:\\
\[  ({\lambda}^1(x),{\lambda}^2(x))\rightarrow ({\xi}(x)=x{\lambda}^1(x)+{\lambda}^2(x),{\xi}_x(x)={\lambda}^1(x),{\xi}_{xx}(x)=0)   \]
The only two non-zero components of the Spencer operator become:\\
\[  {\partial}_x{\xi}(x)-{\xi}_x(x)=x{\partial}_x{\lambda}^1(x)+{\partial}_x{\lambda}^2(x), {\partial}_x{\xi}_x(x)-0={\partial}_x{\xi}_x(x)={\partial}_x{\lambda}^1(x) \]
Equating to zero these two components amounts to have:\\
\[    {\partial}_x{\lambda}^1=0 \hspace{1cm},\hspace{1cm}  {\partial}_x{\lambda}^2=0  \]
Accordingly, gauging $\lambda$ just amounts to choose an arbitrary section of $R_2$.\\
The final touch, that could not be in the mind of any reader even on this very simple example, is to work out the formal adjoint of the Spencer operator. For this, multiplying the first component by a test function ${\sigma}(x)$, the second by a test function ${\mu}(x)$, then summing and integrating by parts, we get the following operator with second members $(f,m)$:\\
\[      {\partial}_x{\sigma}=f  \hspace{1cm},\hspace{1cm}  {\partial}_x\mu +\sigma=m   \]
The comparison with the previous mechanical results needs no comment.\\

Taking into account this example, we now study the foundation of elasticity theory and we restrict the study to 2-dimensional (infinitesimal) elasticity for simplicity as the general situation has already been treated elsewhere and we just want to explain why {\it the only founding problem of elasticity is the choice of an underlying Lie pseudogroup and an adapted differential sequence}.\\

\noindent
1){\it The gauge sequence cannot be used}: \\
Looking at the book [5] written by E. and F. Cosserat, it seems {\it at first sight} that they just construct the first operator of the gauge sequence for $n=1$ ([5], p 7), $n=2$ ([5], p 66), $n=3$ ([5], p 123) and finally $n=4$ ([5], p 189) in the linearized framework. This is {\it not true} indeed because, according to the comment previously done, the adjoint operator is a divergence like operator, a situation not met in the couple-stress equations which is a linear operator with constant coefficients, but not of divergence type. In fact, a carefull study of the book proves that {\it somewhere the action of the group on the space is used}, but this is well hidden among many very technical formulas (Compare [5], p 136 with [21] p 295).\\

\noindent
2){\it The Janet sequence cannot be used}:\\
This result is even more striking because {\it all} texbooks of elasticity use it along the same scheme that we now describe. Indeed, after gauging the translation by defining the "{\it displacement vector}"  $\xi=({\xi}^1(x),{\xi}^2(x))$ of the body, from the initial point $x=(x^1,x^2)$ to the point $y=x+\xi(x)$, one introduces the (small) "{\it deformation tensor}" $\epsilon=1/2{\cal{L}}(\xi)\omega$ as one half the Lie derivative with respect to $\xi$ of the euclidean metric $\omega$, namely, in our case, the three components only (care):\\
\[ \epsilon=({\epsilon}_{11}={\partial}_1{\xi}^1,{\epsilon}_{12}={\epsilon}_{21}=1/2({\partial}_1{\xi}^2+{\partial}_2{\xi}^1),{\epsilon}_{22}={\partial}_2{\xi}^2)\]
From the mathematical point of view, one uses to consider the Lie operator ${\cal{D}}\xi={\cal{L}}(\xi)\omega:T\rightarrow S_2T^*$ (symmetric tensors), sometimes called Killing operator, through the formula:\\
\[      ({\cal{D}}\xi)_{ij}\equiv {\omega}_{rj}{\partial}_i{\xi}^r+{\omega}_{ir}{\partial}_j{\xi}^r+{\xi}^r{\partial}_r{\omega}_{ij}={\Omega}_{ij}=2{\epsilon}_{ij}   \]
One may check at once the only generating "{\it compatibility condition}" ${\cal{D}}_1\epsilon=0$, namely:\\
\[    {\partial}_{11}{\epsilon}_{22}+{\partial}_{22}{\epsilon}_{11}-2{\partial}_{12}{\epsilon}_{12}=0 \]
which is nothing else than the Riemann tensor of a metric, linearized at $\omega$.\\
However, the main experimental reason for introducing the first operator of this type of Janet sequence is the fact that the deformation is made from the displacement and first derivatives but must be invariant under any rigid motion. In the general case it must therefore have $(n+n^2)-(n+n(n-1)/2)=n(n+1)/2$ components, that is 3 when $n=2$, and this is the reason why introducing the deformation tensor $\epsilon$. For most finite element computations, the action density (local free energy) $w$ is a (in general quadratic) function of $\epsilon$ and people use to define the stress by the formula ${\sigma}^{ij}=\partial w/\partial {\epsilon}_{ij}$ which is not correct because $w$ only depends on ${\epsilon}_{11},{\epsilon}_{12}, {\epsilon}_{22}$ when $n=2$  as the deformation tensor is symmetric {\it by construction}. Finally, textbooks escape from this trouble by {\it deciding} that the stress should be symmetric and this is a {\it vicious circle} because we have proved it was not an assumption but an experimental result depending on specific constitutive laws. Accordingly, when $n=2$, we should have ${\sigma}^{ij}{\epsilon}_{ij} = {\sigma}^{11}{\epsilon}_{11}+(2{\sigma}^{12}){\epsilon}_{12}+{\sigma}^{22}{\epsilon}_{22}$. Hence, even if we find the correct stress equations with {\it this convenient duality} keeping the factor "2", we have no way to get the stress {\it and} couple-stress equations {\it together}.\\

\noindent
3){\it Only the Spencer sequence can be used}:\\
Let us construct the formal adjoint of the Spencer operator by multiplying all the $(2\times 2)+2=6$ linearly independent nonzero components by corresponding test functions. For simplifying the summation, we shall raise and lower the indices by means of the (constant) euclidean metric, setting in particular ${\xi}_i={\omega}_{ir}{\xi}^r$ and ${\xi}_{i,j}={\omega}_{ir}{\xi}^r_j$. The only nonzero first jets coming from the $2\times 2$ skewsymmetric infinitesimal rotation matrix of first jets are now ${\xi}_{1,2}=-{\xi}_{2,1}$ while the second order jets are zero because isometries are linear transformations. We obtain in the present situation:\\
\[ {\sigma}^{11}{\partial}_1{\xi}_1+{\sigma}^{12}({\partial}_1{\xi}_2-{\xi}_{1,2})+{\sigma}^{21}({\partial}_2{\xi}_1-{\xi}_{2,1})+{\sigma}^{22}{\partial}_2{\xi}_2+{\mu}^{r,12}{\partial}_r{\xi}_{1,2}  \]
Integrating by parts and changing the sign, we just need to look at the coefficients of ${\xi}_1,{\xi}_2$ and ${\xi}_{1,2}$, namely:\\
\[  \begin{array}{lcl}
 {\xi}_1 & \longrightarrow  &  {\partial}_1{\sigma}^{11}+{\partial}_2{\sigma}^{21}=f^1  \\
 {\xi}_2 & \longrightarrow   & {\partial}_1{\sigma}^{12}+{\partial}_2{\sigma}^{22}=f^2  \\
 {\xi}_{1,2} & \longrightarrow  & {\partial}_r{\mu}^{r,12}+{\sigma}^{12}-{\sigma}^{21}=m^{12}  
\end{array}  \]
in order to get the adjoint operator $ad(D):{\wedge}^{n-1}T^*\otimes R_1^*\rightarrow {\wedge}^nT^*\otimes R_1^*:(\sigma,\mu)\rightarrow (f,m)$ relating for the first time the torsor framework to the Lie coalgebroid $R_1^*$ (see the beginning of section 3). These equations are {\it exactly} the three stress and couple-stress equations of 2-dimensional Cosserat elasticity. In the n-dimensional case, a similar calculation, left to the reader as an exercise of indices, should produce {\it exactly} the $n(n+1)/2$ stress and couple-stress equations in general. It is now possible to enlarge the group in order to get more equations, that is {\it as many equations as the number of group parameters}. Using the conformal group of space-time, the 4 elations give rise to 4 nonzero second order jets only which allow to exhibit the 4  Maxwell equations for the induction $(\vec{H},\vec{D})$ along lines only sketched by H. Weyl in [20] because the needed mathematics were not available before 1970. But, as we already said, this is another story !.\\

\noindent
{\bf REMARK 3.8}: It becomes now clear that the $dim({\wedge}^2T^*\otimes{ \cal{G}})=(n(n-1)/2)(n(n+1)/2)=n^2(n^2-1)/4$ {\it first order} compatibility conditions for the Cosserat fields [3,4,5,11,27] (the so-called torsion and curvature of E. Cartan [2]) are described by the second ({\it first order}) Spencer operator in the Spencer sequence while the $n^2(n^2-1)/12$ {\it second order} compatibility conditions for the deformation tensor (the so-called Riemann curvature) are described by the second operator ${\cal{D}}_1$ in the Janet sequence (see [22], Example 10, p 249-258). Accordingly, {\it the torsion}+{\it curvature of Cartan is not at all the generalization of the curvature of Riemann}, contrary to what is still claimed in mathematical physics today.\\

\noindent
4)  {\bf PARAMETRIZATION}:\\

    The main tool in this last section will be {\it duality theory}, namely the systematic use of the {\it formal adjoint} of an operator. For this, if $E$ is a vector bundle, we introduce its dual $E^*$ to be the vector bundle with inverse transition matrix (for example $T^*$ is the dual of $T$). The formal adjoint of an operator ${\cal{D}}:E\rightarrow F$ is the operator ${ad(\cal{D})}:{\wedge}^nT^*\otimes F^*\rightarrow {\wedge}^nT^*\otimes E^*$ defined by the relation:\\
\[    <\lambda,{\cal{D}}\xi>=<{ad(\cal{D})}\lambda,\xi> + d\omega     \]
where $\lambda$ is a test vector density and $\omega\in {\wedge} ^{n-1}T^*$ comes from Stokes formula of integration by part. Any operator can be considered as the formal adjoint of another operator because we have the identity $ad(ad(\cal{D}))=\cal{D}$. Also, if $P, Q\in D$ and $P=a^{\mu}d_{\mu}$, then $ad(P)=\sum (-1)^{\mid \mu \mid}d_{\mu}a^{\mu}$ and $ad(PQ)=ad(Q)ad(P)$.\\

\noindent
{\bf DEFINITION 4.1}: We now define $N$ from $ad(\cal{D})$ exactly as we already defined $M$ from $\cal{D}$, that is $N=D\lambda/Dad(\cal{D})\lambda$. Of course $N$ {\it highly depends on the presentation of} $M$.\\

   The following nontrivial theorem, first obtained in [18], provides a {\it purely formal test} for deciding about the existence of a parametrization and exhibiting one. At the same time it establishes a link with the two previous sections and is already implemented on the computer algebra package 
http://wwwb.math.rwth-aachen.de/OreModules.\\

\noindent
{\bf THEOREM 4.2}: A test for checking that a given differential module $M$ is torsion-free proceeds in 5 steps:\\
1) Write the defining system as the kernel of an operator $\cal{D}$.\\
2) Construct its formal adjoint $ad(\cal{D})$.\\
3) Work out generating CC for $ad(\cal{D})$ as an operator $ad({\cal{D}}_{-1})$.\\
4) Construct $ad(ad({\cal{D}}_{-1}))={\cal{D}}_{-1}$.\\
5) Work out generating CC for ${\cal{D}}_{-1}$ as an operator ${\cal{D}}'$.\\  
Then $M$ is torsion-free if and only if $\cal{D}$ and ${\cal{D}}'$ have the same solutions (both provide $M$).\\

\noindent
{\bf REMARK 4.3}: We have ${ad({\cal{D}}\circ {\cal{D}}_{-1})=ad({\cal{D}}_{-1})\circ ad(\cal{D})} \equiv 0$ and thus ${\cal{D}}\circ {\cal{D}}_{-1}\equiv 0$. Accordingly $\cal{D}$ ranges among the CC of ${\cal{D}}_{-1}$, a result symbolically written as $\cal{D}\leq {\cal{D}}'$, and $M$ is torsion-free if and only if $\cal{D}={\cal{D}}'$. In that case, with a slight abuse of language, the kernel of $\cal{D}$ is parametrized by the image of ${\cal{D}}_{-1}$. Otherwise, any CC in ${\cal{D}}'$ but not in $\cal{D}$ provides a torsion element of $M$ and ${\cal{D}}_{-1}$ provides a parametrization of the system determined by $M/t(M)$ or, equivalently, $M'=M/t(M)$ is the torsion-free module determined by ${\cal{D}}'$.\\

For the reader aware of homological algebra, the following rather "{\it magic}" corollary depends on difficult technical though classical results of homological algebra ([1] but the best "basic" references are [7], Lemma 3.8, p 147; [16], Theorem 9, p 133; [25], Corollary 6.18, p 186).\\

\noindent
{\bf COROLLARY 4.4}: The differential module $ext^1_D(N,D)=t(M)$ does not depend on the presentation of $N$.\\

\noindent
{\bf EXAMPLE 4.5}: As a first striking result, that does not seem to have been noticed by mechanicians up till now, let us consider the situation of classical elasticity theory where $\cal{D}$ is the Killing operator for the euclidean metric, namely $\cal{D}\xi=\cal{L}(\xi)\omega$ and ${\cal{D}}_1$ the corresponding CC, namely the linearized Riemann curvature with $n^2(n^2-1)/12$ components. In that case, as it is well known that the Poincar\'{e} sequence for the exterior derivative is self-adjoint {\it up to sign} (for $n=3$ the adjoints of $grad,curl,div$ are respectively $div,curl,grad$) then {\it the first extension module does not depend on the differential sequence used} and therefore vanishes. Accordingly, $ad(\cal{D})$ generates the CC of $ad({\cal{D}}_1)$. Hence, in order to parametrize the stress equations, that is $ad(\cal{D})$, one just needs to compute $ad({\cal{D}}_1)$. For $n=2$, we get:\\
\[    \phi ({\partial}_{11}{\epsilon}_{22}+{\partial}_{22}{\epsilon}_{11}-2{\partial}_{12}{\epsilon}_{12})={\partial}_{22}\phi{\epsilon}_{11}-2{\partial}_{12}\phi{\epsilon}_{12}+{\partial}_{11}\phi{\epsilon}_{22}+{\partial}_1(...)+{\partial}_2(...)   \]
and recover ... the parametrization by means of the Airy function in a rather unexpected way.\\
   Exhibiting a parametrization for $n\geq 3$ thus becomes a straightforward exercise of computer algebra, the number of (pseudo)-potentials being $n^2(n^2-1)/12$ along the last remark of the preceding section.\\
   
{\bf EXAMPLE 4.6}: Contrary to the preceding situation, the second order Einstein operator $\cal{D}$ in General Relativity has $n/(n+1)/2$ components and is self-adjoint as it comes from Hilbert variational calculus of the scalar Riemann density. The corresponding CC are the well known $n$ divergence conditions and the adjoint ${\cal{D}}_{-1}$ is therefore the Killing operator for the Minkowski metric. According to the previous example, its CC ${\cal{D}}'$ is thus the linearized Riemann curvature with $n^2(n^2-1)/12$ components. It follows that $\cal{D}\leq {\cal{D}}'$ {\it strictly} when $n\geq 4$ and no parametrization can be found in that case.\\

\noindent
{\bf EXAMPLE 4.7}: We finally treat the case of Cosserat equations. For this, instead of using the Janet sequence as before, we now use the Spencer sequence which is isomorphic to the gauge sequence though with quite different operators. However, according to the general theorems of homological algebra, the existence of a parametrization does not depend on the differential sequence used and therefore follows again, like in the first example of this section, from the fact that the Poincar\' {e} sequence is self-adjoint up to the sign. In the present situation, we have $C_r={\wedge}^rT^*\otimes R_1\simeq {\wedge}^rT^*\otimes \cal{G}$ with ${dim(\cal{G})}=n(n+1)/2$. As we have shown in the last section that the Cosserat equations were just $ad(D_1)$, their {\it first order} parametrization is thus described by $ad(D_2)$ and needs $dim(C_2)=n^2(n^2-1)/4$ (pseudo)-potentials according to the last remark of the preceding section. We provide the details when $n=2$ but we know at once that we must use $3$ (pseudo)-potentials only. The case $n=3$ could be treated similarly and is left to the reader as an exercise.\\
   The Spencer operator $D_1$ is described by the equations:\\
\[  {\partial}_1{\xi}_1=A_{11}, {\partial}_1{\xi}_2-{\xi}_{1,2}=A_{12}, {\partial}_2{\xi}_1-{\xi}_{2,1}=A_{21}, {\partial}_2{\xi}_2=A_{22}, {\partial}_1{\xi}_{1,2}=B_1, {\partial}_2{\xi}_{1,2}=B_2  \]
because $R_1$ is defined by the equations  ${\xi}_{1,1}=0,{\xi}_{1,2}+{\xi}_{2,1}=0, {\xi}_{2,2}=0$. \\
   Accordingly the 3 CC describing the Spencer operator $D_2$ are:\\
\[  {\partial}_2A_{11}-{\partial}_1A_{21}+B_1=0, {\partial}_2A_{12}-{\partial}_1A_{22}+B_2=0, {\partial}_2B_1-{\partial}_1B_2=0  \]
Multiplying these equations respectively by ${\phi}^1, {\phi}^2, {\phi}^3$, then summing and integrating by part, we get $ad(D_2)$ and the desired first order parametrization in the form:\\
\[   {\sigma}^{1,1}=-{\partial}_2{\phi}^1, {\sigma}^{1,2}=-{\partial}_2{\phi}^2, {\sigma}^{2,1}={\partial}_1{\phi}^1, {\sigma}^{2,2}={\partial}_1{\phi}^2, {\mu}^{1,12}=-{\partial}_2{\phi}^3+{\phi}^1, {\mu}^{2,12}={\partial}_1{\phi}^3+{\phi}^2  \]
as announced previously. It is important to notice that such a parametrization, which could also be obtined by localization, is coherent with the classical one already obtained by localization in Example 2.6 and which can be recovered if we cancel the couple-stress.\\

\noindent
5)  {\bf CONCLUSION}:\\

    We have already proved in our book "{\it Lie Pseudogroups and Mechanics}" (1988) that the group foundation of elasticity pioneered by E. and F. Cosserat (1909) was just described by adding to the {\it non-linear Spencer sequence} (1972) a convenient variational calculus. As a crucial conclusion and similarly to classical elasticity, even if the initial framework is non-linear, the resulting stress/couple-stress equations are linear. Accordingly, in order to parametrize these equations, one just needs to refer to {\it infinitesimal Lie equations} and the three corresponding {\it canonical linear differential sequences} that can be constructed, namely the {\it gauge sequence}, the {\it Janet sequence} and the {\it Spencer sequence}, though only the last one is useful.\\
    The main result of this paper, coming from {\it unavoidable arguments} of homological algebra, is that the Cosserat equations are just described by the formal adjoint of the first Spencer operator while their parametrization is provided by the formal adjoint of the second Spencer operator. This result finally explains why the parametrization of the stress/couple-stress equations is {\it first order} while the parametrization of the classical stress is {\it second order} as it comes from dualizing the Janet sequence.\\
     We are happy to to pay such a tribute to E. and F. Cosserat on the occasion of this anniversary and hope it will give a new impulse for generalizing their work.\\

\noindent
{\bf BIBLIOGRAPHY}:\\

\noindent
[1]  BOURBAKI: Alg\`{e}bre, Chapitre 10, Alg\`{e}bre Homologique, Masson, Paris, 1980.\\
\noindent
[2] CARTAN, E.: Sur une g\' {e}n\' {e}ralisation de la notion de courbure de Riemann et les espaces \`{a} torsion, C. R. Acad\' {e}mie des Sciences Paris, 174, 1922, P. 522.\\
\noindent
[3] CHWOLSON,O.D.: Trait\' {e} de physique.(In particular III,2,p. 537 + III,3, p. 994 + V, p. 209), Hermann, Paris, 1914.\\
\noindent
[4] COSSERAT, E. and F.: Note sur la th\' {e}orie de l'action euclidienne. In : Trait\' {e} de m\' {e}canique rationelle (Appell, P., ed.),t. III, pp. 557-629, Gauthiers-Villars, Paris, 1909.\\
\noindent
[5]  COSSERAT, E. and F.: Th\'{e}orie des Corps D\'{e}formables, Hermann, Paris, 1909.\\
\noindent
[6]  GR\"{O}BNER, W.: \"{U}ber die Algebraischen Eigenschaften der Integrale 
von Linearen Differentialgleichungen mit Konstanten Koeffizienten, 
Monatsh. der Math., 47, 1939, 247-284.\\
\noindent
[7]  HU, S.T.: Introduction to Homological Algebra, Holden-Day, 1968.\\
\noindent
[8]  JANET, M.: Sur les Syst\`{e}mes aux d\'{e}riv\'{e}es partielles, 
Journal de Math., 8, 3, 1920, 65-151.\\
\noindent
[9]  KALMAN, E.R., YO, Y.C., NARENDA, K.S.: Controllability of Linear Dynamical
Systems, Contrib. Diff. Equations, 1, 2, 1963, 189-213.\\
\noindent
[10]  KASHIWARA, M.: Algebraic Study of Systems of Partial Differential 
Equations, M\'emoires de la Soci\'et\'e Math\'ematique de France 63, 1995, 
(Transl. from Japanese of his 1970 Master's Thesis).\\
\noindent
[11] KOENIG, G.: Le\c{c}ons de cin\' ematique (The Note "Sur la cin\' {e}matique d'un milieu continu" by E. and F. Cosserat, pp. 391-417, has never been quoted elsewhere), Hermann, Paris, 1897.\\
\noindent
[12] KUMPERA, A., SPENCER, D. C.: Lie equations, ANN. Math. Studies 73, Princeton University 
Press, Princeton, 1972.\\
\noindent
[13]  KUNZ, E.: Introduction to Commutative Algebra and Algebraic Geometry, 
Birkh\"{a}user, 1985.\\
\noindent
[14]  MAISONOBE, P., SABBAH, C.: D-Modules Coh\'erents et Holonomes, Travaux en
Cours, 45, Hermann, Paris, 1993.\\
\noindent
[15]  MALGRANGE, B.: Cohomologie de Spencer (d'apr\`{e}s Quillen), 
S\'em. Math. Orsay, 1966.\\
\noindent
[16]  NORTHCOTT, D.G.: An Introduction to Homological Algebra, Cambridge University Press, 1966.\\
\noindent
[17]  PALAMODOV, V.P.: Linear Differential Operators with Constant Coefficients,
Grundlehren der Mathematischen Wissenschaften 168, Springer, 1970.\\
\noindent
[18]  POMMARET, J.-F.: Dualit\'{e} Diff\'{e}rentielle et Applications,
 C. R. Acad. Sci. Paris, 320, S\'erie I, 1995, 1225-1230.\\
\noindent
[19]  POMMARET, J.-F.: Fran\c{c}ois Cosserat et le Secret de la Th\'{e}orie Math\'{e}matique de 
l'Elasticit\'{e}, Annales des Ponts et Chauss\'{e}es, 82, 1997, 59-66.\\
\noindent
[20] POMMARET, J.-F.: Group interpretation of coupling phenomena, Acta Mechanica, 149, 2001,pp. 23-39.\\
\noindent
[21]  POMMARET, J.-F.: Lie Pseudogroups and Mechanics, Gordon and Breach, New York, 1988.\\
\noindent
[22]  POMMARET, J.-F.: Partial Differential Equations and Group Theory: New
Perspectives for Applications, Kluwer, 1994.\\
\noindent
[23]  POMMARET, J.-F.: Partial Differential Control Theory, Kluwer, 2001.\\
(http://cermics.enpc.fr/$\sim$pommaret/home.html)\\
\noindent
[24]  POMMARET, J.-F.Algebraic Analysis of Control Systems Defined by Partial Differential Equations, in Advanced Topics in Control Systems Theory, Lecture Notes in Control and Information Sciences 311, Chapter 5, Springer, 2005, 155-223.\\
\noindent
[25]  ROTMAN, J.J.: An Introduction to Homological Algebra, Pure and Applied 
Mathematics, Academic Press, 1979.\\
\noindent
[26]  SPENCER, D.C.: Overdetermined Systems of Partial Differential Equations, 
Bull. Amer. Math. Soc., 75, 1965, 1-114.\\
\noindent
[27] TEODORESCU, P. P. : Dynamics of linear elastic bodies, Editura Academiei, Bucuresti, Romania; Abacus Press, Tunbridge, Wells, 1975.\\
\noindent
[28] WEYL, H. : Space, Time, Matter, Springer, Berlin, 1918, 1958; Dover, 1952.\\
\noindent
[29] YANG, C. N. : Magnetic monopoles, fiber bundles and gauge fields, Ann. New York Acad. Sciences, 294, 1977, p. 86.\\
\noindent
[30]  ZERZ, E.:Topics in Multidimensional Linear Systems Theory, Lecture 
Notes in Control and Information Sciences 256, Springer, 2000.\\

\end{document}